\documentclass[11pt,reqno]{amsart}

\usepackage{amsmath,amsfonts,amssymb,amsthm,amscd}

\makeindex

\usepackage{amssymb}
\usepackage[latin1]{inputenc}
\usepackage{ epsfig,epsf}
\usepackage{enumerate}
\usepackage{indentfirst}
\usepackage{euscript}
\usepackage{graphicx}
\usepackage{amsmath}
\usepackage{amsfonts}
\usepackage{amssymb}
\makeindex

 \theoremstyle{plain}

\newtheorem{proposition}{Proposition}

\newtheorem{remark}{Remark}
\theoremstyle{definition}

\newtheorem{theorem}{Theorem}

\usepackage{amsmath,amsfonts,  amssymb,amsthm,amscd}
\makeindex
\usepackage{amssymb}
\usepackage{graphics}

\usepackage[latin1]{inputenc}
\usepackage{epsf}

\usepackage{ epsfig,epsf}
\usepackage{enumerate}
\usepackage{indentfirst}
\usepackage{euscript}
\usepackage{graphicx}
\usepackage{amsmath}
\usepackage{amsfonts}
\usepackage{amssymb}
\makeindex

 \theoremstyle{plain}

 \title[ Codimension two Umbilic points  ]{ Codimension two Umbilic points
  on  Surfaces  Immersed  in $\mathbb R^3$}

 \author[R. Garcia]{Ronaldo Garcia}

 \author[J. Sotomayor]{Jorge Sotomayor}


 \thanks{Both  authors are fellows of  CNPq.
 This work was done under the project CNPq PADCT 620029/2004-8 and CNPq
 473824/04-3.
The first author was partially supported by  FUNAPE/UFG}

\begin{document}
 \maketitle

 \begin{abstract} In this paper is studied the behavior of lines of curvature
near  umbilic points that appear generically on surfaces depending
on two parameters.
  \end{abstract}

\section{ Introduction } \label{sec:1}

The study of principal curvature lines, along which surfaces
immersed in  ${\mathbb R}^3$ bend maximally and minimally,  goes
back to the work of Monge \cite {Mo}. There was determined the
behavior of principal lines  on  the  triaxial ellipsoids. The
analysis around  umbilic points, at which surfaces bend equally in
all directions,
 had a preponderant  role in this work, which may be
  regarded as   the first one on the
 subject of foliations with singularities.

\begin{figure}[htbp]
  \begin{center}
  \hskip .5cm
 \includegraphics[angle=0, width=10cm]{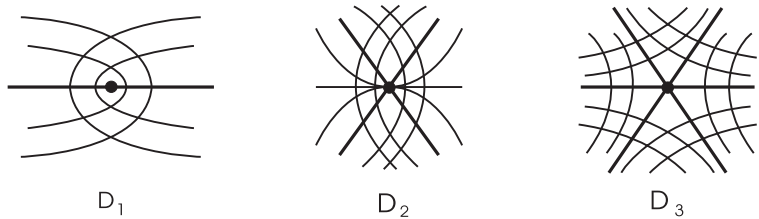}
  \caption{Principal curvature lines near the
Darbouxian $D_i ,  \, i\,=\,1, \,2, \,3,$
 umbilic points  and their separatrices \label{fig:darboux}}
 \end{center}
  \end{figure}

In 1896 Darboux \cite{da} carried out
a complete local analysis of  umbilic points in the class of analytic
surfaces, under generic conditions on the third derivatives  of the
immersions of surfaces at the umbilic points.
 Darboux established that there are three
 generic patterns,    called here {\it Darbouxian Umbilics }
$D_1$, $D_2$ and $D_3$, illustrated in Fig.
 \ref{fig:darboux}. This figure represents both families of
 principal curvature lines, minimal $\mathcal F_1$ and  maximal $\mathcal F_2$, associated to minimal
 and maximal principal curvatures $k_1$ and $k_2$,  $k_1< k_2$, expect at umbilic points.
The subscript stands for the number of umbilic separatrices
approaching the point. In this case the number is the same for
both the minimal and maximal
 principal curvature foliations.

This result of Darboux was rediscovered by   Gutierrez and
Sotomayor, \cite{gs1, gs2} in the context  of Structural Stability
of principal lines on surfaces of class $C^4$, $r\geq 4$.  They
showed that the Darbouxian
 umbilic points are generic, corroborating in $C^4$, $r\geq 4$, the above mentioned
 results of Darboux,   and characterize the structurally stable ones, under  small $C^3$
 deformations of the surface. See also the work of Bruce and Fidal \cite{bf}.
 In Sotomayor  \cite{sho}  can be found comments on the early developments on
lines of curvature and  their relation with recent works on
foliations with singularities geometrically defined   on surfaces
immersed in ${\mathbb R}^3$.

 In \cite{ggs, gs4} the simplest bifurcation patterns of Darbouxian
umbilics were established by Garcia, Gutierrez and Sotomayor. This
corresponds to the weakening of the Darbouxian structural
stability  conditions in the most stable and generic way,
 leading to {\it codimension one umbilic points} and their generic
bifurcations.
See section \ref {sec:2} for a review of the analytic expressions of these conditions.

\begin{figure}[htbp]
  \begin{center}
  \hskip 1cm
 \includegraphics[angle=0, width=7cm]{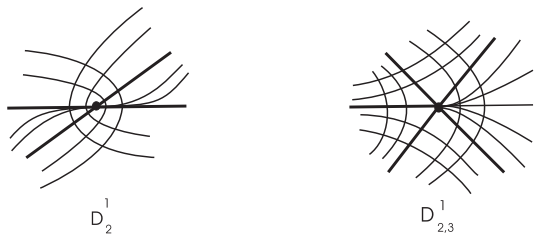}
 \caption{Lines of Curvature near Codimension one Umbilic Points $D^{1}_{2}$, left,
   and  $D^{1}_{2,3}$, right.  \label{fig:semid}}
    \end{center}
  \end{figure}

The expression {\it codimension one} means that the umbilic points
appear generically on  one-parameter families of immersed
surfaces. This corresponds to  families  in general position with
respect
 to a regular part of   the bifurcation set --a hypersurface-- in the
 infinite dimensional space of
 immersions. The two patterns  of umbilics of
codimension one established in  \cite{ggs, gs4},
 denominated $D^{1}_{2}$ and $D^{1}_{2,3}$, are illustrated in Fig.
\ref{fig:semid}.  The superscript stands for the codimension and the subscripts stand
for the number
 of separatrices approaching the umbilic point. In the first case,
 the number of separatrices is the same
 for both the minimal and maximal principal curvature foliations.
 In the second case, they are not
 equal and, in our notation,  appear  separated by a comma.
 See Theorem \ref{th:2} in section \ref{sec:2}.

 This paper pursues the
 analysis in \cite{ggs}, increasing in one  unity
the number of parameters on which the immersion of the surface
depends. The general {\it codimension two} umbilic points
exhibited by  generic bi-parametric families of surfaces immersed
in ${\mathbb R}^3$, consisting on four patterns  -- called
$D^{2}_{1}\,$, $D^{2}_{2p} \,$, $D^{2}_{3}\,$ and $D^{2}_{2h}$
 {\it codimension two} umbilic points --  will be established
here. The superscript stands for the codimension. The symbols $p$,
for {\it parabolic},  and $h$, for {\it hyperbolic}, have been
added to the subscripts above in order to distinguish types that
are not discriminated only by the number of separatrices.

The first three patterns
 can be described  as follows:
$D^{2}_{1}$ is
 topologically equivalent to a point    $D_{1}$;
$D^{2}_{2p}$ is topologically equivalent to a point
 $D_{2}$; $D^{2}_{3}$ is topologically equivalent to  a point
 $D_{3}$. See Fig. \ref{fig:darboux}.

\begin{figure}[htbp]
  \begin{center}
  \hskip 1cm
 \includegraphics[angle=0, width=7cm]{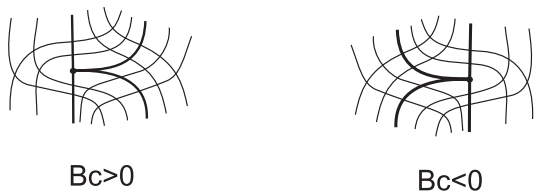}
  \caption{Lines of curvature near the umbilic
   point $D^{2}_{2h}$ \label{fig:d13}}
    \end{center}
  \end{figure}

The pattern $D^{2}_{2h}$ is illustrated in Fig. \ref{fig:d13},
where $c$ and $B$ are coefficients  defined in equation
(\ref{eq:1}), pertinent to a Monge presentation of the immersion
at the umbilic point.  This pattern is new with respect to
previous types. Also no proof of its structure  has been found in
the literature.

Section \ref{sec:2} is devoted to review some  background on
umbilic points of codimension less than two. The local
configurations of principal lines at umbilic points of types
$D^{2}_{1}\,$, $D^{2}_{2p} \,$, $D^{2}_{3}\,$ and $D^{2}_{2h}$
  are established in Section
\ref{sec:3}. The fact that these points actually are those of {\it
codimension two}, in a geometric sense, is established in Section
\ref{sec:4}.

 The study of the
bifurcation diagrams of lines of principal curvature lines at
umbilic points for bi-parametric families of surfaces in ${\mathbb
R}^3$, roughly prepared in Section \ref{sec:4},  involve some
technical and lengthly work, especially in the case of
$D^{2}_{2h}$, and  will be postponed to a future paper
\cite{bif2}.

\section{ Preliminaries on   Umbilic Points} \label{sec:2}
The following assumptions will hold from now on.

 Let $p_0$ be an
umbilic point of an immersion $\alpha$ of an oriented surface
$\mathbb M $ into $\mathbb R ^3$, with a once for all fixed
orientation. It will be assumed that $\alpha$ is of class $C^k ,\;
k \geq \, 6$.
 In a local Monge  chart
 near $p_0$  and a positive  adapted $3-$frame,
  $\alpha$ is given by $\alpha (u,v)=(u,v,h(u,v))$,
where

 \begin{equation}\label{eq:1}
  \aligned h (u,v) &= \frac k2 (u^2+v^2) + \frac a6 u^3 +\frac
b2 u v^2 +\frac c6 v^3
+\frac A{24} u^4 + \frac B6 u^3 v\\
 &+\frac C4 u^2v^2 + \frac D6 u v^3 + \frac E{24} v^2
+\frac{a_{50}}{120}u^5+ \frac{a_{41}}{24}u^4v\\
&+ \frac{a_{32}}{12}u^3v^2 +\frac{a_{23}}{12}u^2v^3+
\frac{a_{14}}{24}uv^4+\frac{a_{05}}{120}v^5+ h.o.t
\endaligned
\end{equation}

Notice that, without loss of generality,  the term $u^2 v$ has
been eliminated from this expression by means of a rotation in the
$(u,v)-$ frame.

 According to \cite{sp} and \cite{st}, the differential equation of lines of
curvature in terms of  $I=Edu^2 +2Fdudv+Gdv^2$ and  $II=edu^2
+2fdudv+gdv^2$ around $p_0$ is:

\begin{equation}\label{eq:lc}
(Fg-Gf)dv^2+(Eg-Ge)dudv+(Ef-Fe)du^2=0.
\end{equation}

\noindent  Therefore the  functions $L=Fg-Gf$, $M=Eg-Ge$ and
$N=Ef-Fe$ are:
$$\aligned L &= h_uh_v h_{vv}- (1+h_v^2)h_{uv}\\
M &= (1+h_u^2)h_{vv}- (1+h_v^2)h_{uu}\\
N &= (1+h_u^2)h_{uv} - h_u h_v h_{uu}.\endaligned $$

Extensive calculation gives

\begin{equation}\label{eq:lmn}
\aligned L &= -b v -\frac 12 B
u^2 - (C-k^3)u v -\frac 12 D v^2
-\frac{a_{41}}6u^3\\
&+\frac 12(4bk^2+k^2a-a_{32})u^2v
+\frac 12(3k^2c-a_{23})uv^2 -\frac 16(a_{14}+3bk^2)v^3
 + h.o.t\\
 M &= (b-a)u + c v
+ \frac 12( C-A+2k^3)u^2 + (D-B) u v+ \frac 12 ({E-C} -2k^3)v^2
\\ & +\frac 16[6bk^2(a+b) +a_{32}-a_{50}]u^3+
\frac 12(a_{23}+2ck^2-a_{41})u^2v\\
&+\frac 12[a_{14}-a_{32}-2 k^2(a+b)]uv^2+
\frac 16(a_{05}-a_{23}-6ck^2)v^3+
 h.o.t\\
N&= b v + \frac 12 B  u^2 + (C-k^3)u v + \frac 12 D  v^2
+ \frac 16  a_{41} u^3\\
&+\frac 12(a_{32}-3ak^2)u^2v
+\frac 12( a_{23}-k^2c)uv^2
+\frac 16(a_{14}-3bk^2)v^3
+ h.o.t \endaligned \end{equation}

These expressions
 appear, up to order two, in Darboux \cite{da},
including also the coefficient of $u^2 v$, annihilated here by a
rotation in the frame $(u,v)$.

\subsection  {Darbouxian Umbilic Points}\label{ssec:dar}

 For the sake of completeness we will state now  the
  result of Gutierrez and Sotomayor
 in \cite{gs1, gs2}  about the local behavior of
curvature lines near a Darbouxian umbilic point.

\begin{theorem} \label{th:1}
Let $p_0$ be an umbilic point and consider
 the Monge chart as in equation (\ref{eq:1}).

 Suppose that the transversality condition:

$T)\;\; b(b-a)\ne 0$

and  one   the following --Darbouxian-- conditions holds.

\begin{itemize}
\item[{$D_1$})] \quad $ (\frac{c}{2b})^2 -\frac ab +2 <0$
\item[{$D_2$)}]\quad $ (\frac{c}{2b})^2+2 >\frac ab  >1, \quad
a\ne 2b$
 \item[{$D_3$})]\quad $\frac ab <1$
\end{itemize}

Then the behavior of lines of curvature near the umbilic point,
separated in two principal foliations $\mathcal F_1$ and $\mathcal
F_2$, with a common singularity at
 $p_0$ is  as illustrated in  Fig.
\ref{fig:darboux}.
\end{theorem}

The local configuration near umbilics is explained in terms of the
phase portraits of the  singularities of the Lie-Cartan vector
field
 suspension, $X_{\mathcal F},$ of the implicit differential
  equation ${\mathcal F}(u,v,p)=0$, where,

\begin{equation}\label{eq:f}
\aligned
{\mathcal F}(u,v,p)&=L p^2+M p+N=0,\; p=\frac {dv}{du}\\
 X_{\mathcal
F}&=({\mathcal F}_p, p {\mathcal F}_p,-({\mathcal F}_u+p {\mathcal
F}_v)). \endaligned
\end {equation}

\noindent

The projection  by $\pi(u,v,p)=(u,v)$ of the integral curves of
$X_{\mathcal F} $, restricted to ${\mathcal F}(u,v,p) = 0$ are the
lines of curvature near the umbilic point. The singularities of
$X_{\mathcal F}$, restricted to the invariant
 regular surface ${\mathcal F}(u,v,p)=0$, are given by $(0,0,p_i)$ where $p_i$ is
  a root of the equation $({\mathcal F}_u+p{\mathcal
F}_v)(0,0,p)\,=\, p(bp^2-cp+a-2b)=0$. In in the present Darbouxian
case, these singularities are  hyperbolic
 saddles and nodes.

\begin{figure}[htbp]
  \begin{center}
  \hskip .5cm
  \includegraphics[angle=0, width=8cm]{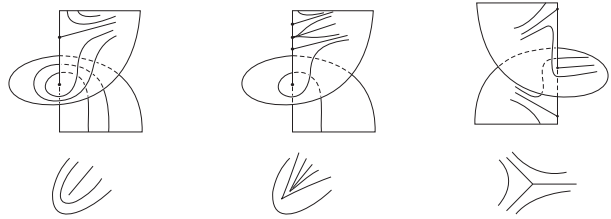}
  \caption{Lie-Cartan  Suspension of  the $D_i$ umbilic points
 \label{fig:darbouxp}}
   \end{center}
   \end{figure}

In  Fig. \ref{fig:darbouxp} is illustrated the behavior of
 $X_{\mathcal F}$ in each of the Darbouxian cases. Here only one
 of the principal foliations has been drawn, for the sake of simplicity.
 The other one is the projection of the other halph cylinder,
 not sketched in the Figure.

\begin{remark} With reference to the illustrations in Fig.
\ref{fig:darboux} of principal curvature foliations, arises the
question of  deciding  which foliation is  minimal and which one
is maximal.

The procedure to follow, going back to the basic definitions, is
to evaluate the normal curvature function
$k_n(p_0,[du:dv])=II(p_0,[du:dv])/I(p_0,[du:dv])$, on the surface
${\mathcal F}(p_0,[du:dv])=0$. Locally, near the projective line
${\mathcal P}=\{(0,0,p)\}$ the set $\mathcal F^{-1}(0)\setminus
\mathcal P$ has two open connected components defined by $
{\mathcal C}_1= \{(p_0, [du:dv]) :
k_n(p_0,[du:dv])=k_1(p_0)\}\setminus \mathcal P $ and ${\mathcal
C}_2= \{( p_0, [du:dv]) : k_n(p_0,[du:dv])=k_2(p_0)\}\setminus
\mathcal P$. The projection of the phase portrait of $X_\mathcal
F$ in the component ${\mathcal C}_1$ (resp. ${\mathcal C}_2$)
defines the minimal (resp. maximal) principal foliation.

This procedure applies to all isolated umbilics and will not be
repeated for those studied below.

\end{remark}

\subsection  {Codimension One Umbilic Points}\label{ssec:sdar}

For future reference, we review a  basic  result of \cite{ggs,
gs4} for  umbilics in generic families
 of surfaces depending on one parameter.

\begin{theorem}  {\cite{ggs, gs4}}\label{th:2}
Let $p_0$ be an umbilic point and consider
$\alpha(u,v)=(u,v,h(u,v))$
 as in equation (\ref{eq:1}).
 Suppose the following conditions hold:

\begin{itemize}
   \item[{$D^{1}_{2}$})]\quad $ cb(b-a)\ne 0\,$  and  either
\quad $ (\frac{c}{2b})^2 - \frac ab +1 = 0\quad$   or \quad $
a=2b.$
 \item[{$D^{1}_{2,3}$})] \quad $b=a\ne 0\quad \text{and} \quad
   \chi={cB}  -({C-A}  + 2k^3)b\ne 0. $
\end{itemize}

 Then the behavior of lines of
curvature near the  umbilic point $p_0$ in  cases $D^{1}_{2}$ and
$D^{1}_{2,3}$,
 is as illustrated  in Fig. \ref{fig:semid}.

\end{theorem}

\begin{figure}[htbp]
  \begin{center}
  \hskip 1cm
  \includegraphics[angle=0, width=12cm]{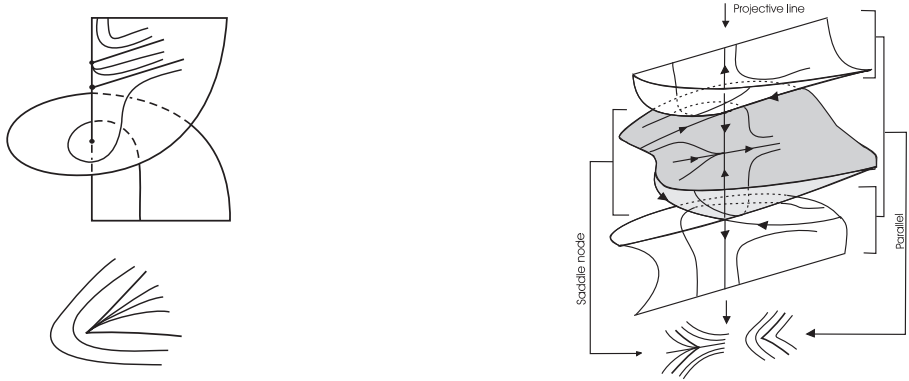}
  \caption{Lie-Cartan suspension  $D^{1}_{2}$, left, and
  $D^{1}_{2,3}$, right.
  \label{fig:lc}}
    \end{center}
  \end{figure}

In  Fig. \ref{fig:lc} is illustrated the behavior of
 $X_{\mathcal F}$ in each of the semi Darbouxian cases.

In this section,  umbilic points which appear generically and
stably on one-parameter families of surfaces have been described
in terms
 of hyperbolic and partially hyperbolic  singularities of vector
fields. The analogy between the structurally stable plane vector
field  singularities  and stable umbilic points on surfaces is
complete in the Darbouxian case. In the  $D^{1}_{2}$  and
$D^{1}_{2,3}$, codimension one  umbilics, the analogy is only
partial.  The full analysis depends both  on  the classic
saddle-node, in the plane,  and  on the phase portrait of a
singularity of a vector field on a conic Morse critical point on
the surface $\mathcal F ^{-1}(0)$.

\section{Patterns of Codimension Two Umbilic Points} \label{sec:3}

In this section will be studied the four simplest  umbilic points,
denominated $D_1^2$, $D^{2}_{2p}$, $D^{2}_{3}$ and $D^{2}_{2h}$,
not discussed previously, assuming that the non-degeneracy
conditions imposed in Theorem \ref{th:2} of section \ref{sec:2},
are violated in the mildest possible fashion and are substituted
by suitable higher order hypotheses. This procedure goes back to
the notions of {\it higher order stability} due to Andronov and
Leontovich \cite{ap}. Se also  \cite{ggs} and \cite{gg}. Here,
referring to equation \ref{eq:1},    $j^3 h(0)$ does not vanish
and the hypotheses are expressed in terms of $j^4 h(0)$ and $j^5
h(0)$.

\subsection  { $D^{2}_{1}$ Umbilic Points}\label{ssec:d21}

\begin{proposition} \label{prop:1} Let $p_0$ be an umbilic point
 and $\alpha(u,v)=(u,v,h(u,v))$
 as in equation \ref{eq:1}.

 If $D^{2}_{1}$): $\;\;c=0$\;\; and\;\; $a=2b\ne 0$,

\noindent then  the configuration  of  principal lines  near $p_0$
is topologically equivalent to    that of  a  Darbouxian $D_1$
 umbilic point. See  Fig. \ref{fig:darboux}, left.
\end{proposition}

\begin{proof}
The differential equation of curvature lines is given by equation
(\ref{eq:f}), where, for the present case:

$$\aligned L&=-bv-\frac 12 Bu^2+(k^3-C)uv-\frac 12 Dv^2 +h.o.t\\
 M &= -bu+\frac 12 (2k^3-A+C)u^2+(D-B)uv+\frac 12(E-2k^3-C) v^2 +h.o.t\\
 N &=bv+\frac 12 Bu^2+(-k^3+C)uv+\frac 12 Dv^2 +h.o.t\endaligned $$

Direct calculation shows that ${\mathcal F}_v(0,0,0)=b\ne 0$. The
solution, $v=v(u,p)$ of the implicit equation
 ${\mathcal
F}(u,v(u,p),p)=0$ near $(0,0,0)$ is given by:

$$v(u,p)=up -\frac B{2b}u^2+  O(3)$$

Near the origin, in the chart $(u,p)$, the vector field
$X_{\mathcal F}$ is given by

 $$\aligned \dot u &= {\mathcal F}_p(u,v(u,p),p) =  -bu+O(2)\\
\dot p &=-({\mathcal F}_u+p{\mathcal F}_v)(u,v(u,p),p)=u[-B+
O(1)]+ bp^3.
\endaligned
 $$

Therefore the point $0$ is a topological saddle of order three,
having the projective line as the unique center manifold.

The analysis of $X_{\mathcal F}$ outside a neighborhood of the
$(0,0,0)$ is trivial, since the three singular points on the
$p-$axis are concentrated at the origin. In the chart $(u,v,q),\;
q=du/dv,$ the correspondent Lie-Cartan vector field is regular at
$(0,0,0)$. Therefore, gluing the phase portraits as in Darbouxian
$D_1$ umbilic point,  the result follows. See Fig.
\ref{fig:darbouxp}, left.
\end{proof}

\subsection  { $D^{2}_{2p}$ and  $D^{2}_{3}$ Umbilic Points}\label{ssec:d23}

\begin{proposition}\label{prop:2} Let $p_0$ be an umbilic point and $\alpha(u,v)=(u,v,h(u,v))$ be
  as in equation \ref{eq:1}.

  If $a=b\ne 0$,
$ \chi = cB- b(C-A +2 k^3)=0 $
 and

\noindent $\xi=  12k^2 b^3+ (a_{32}-a_{50}  )b^2+( 3B^2-3BD
-ca_{41})b+ 3 cB(C-k^3)  \ne 0, \, $

\noindent then the principal configurations of lines of curvature
fall into one of the  two cases:
\begin{itemize}
\item[i)]
 Case $ D^{2}_{2p}$: $\xi b<0 $, which is topologically a $D_2$  umbilic and
\item[ii)] Case $D^{2}_{3}$: $\xi b>0$, which is topologically a
 $D_3$  umbilic.
\end{itemize}
See Fig. \ref{fig:darboux},  center and right , respectively.
\end{proposition}

\begin{proof}
The differential equation of curvature lines is given by
${\mathcal F}(u,v,p)=Lp^2+Mp+N=0,\; p=\frac {dv}{du}$
 where the coefficients are given from equation \ref{eq:lmn} as follows:

\begin{equation}\label{eq:lmn1}
\aligned L &= -b v -\frac 12 B
u^2 - (C-k^3)u v -\frac 12 D v^2 -\frac 16 a_{41} u^3+\frac 12(5bk^2 - a_{32})u^2v\\
&+\frac 12(3k^2c-a_{23})uv^2 -\frac 16(a_{14}+3bk^2)v^3
 + h.o.t\\
 M &= c v
+ \frac 12( C-A+2k^3)u^2 + (D-B) u v+ \frac 12 ({E-C} -2k^3)v^2
\\ & +\frac 16[12bk^2  +a_{32}-a_{50}]u^3+
\frac 12(a_{23}+2ck^2-a_{41})u^2v\\
&+\frac 12[a_{14}-a_{32}-4 k^2b]uv^2+ \frac
16(a_{05}-a_{23}-6ck^2)v^3+
 h.o.t\\
N&= b v + \frac 12 B  u^2 + (C-k^3)u v + \frac 12 D  v^2
+\frac 16  a_{41} u^3+\frac 12(a_{32}-3bk^2)u^2v\\
&+\frac 12( a_{23}-k^2c)uv^2
+\frac 16(a_{14}-3bk^2)v^3
+ h.o.t \endaligned \end{equation}

Direct calculation shows that ${\mathcal F}_v(0,0,0)=b\ne 0$. The
solution  $v=v(u,p)$ of  equation ${\mathcal F}(u,v(u,p),p)=0$ is
given by:

$$v(u,p)=-\frac B{2b}u^2+[\frac{ 3 B (C-k^3)-a_{41}b }{6b^2}]{u^3}+
[ \frac{cB-( C-A+2k^3)b  }{2b^2}]{u^2p}+ O(4)$$

Therefore, in the chart $(u,p)$, the Lie-Cartan vector field $X_{\mathcal F}$ is given by

$$\aligned \dot u &= {\mathcal F}_p(u,v(u,p),p)\\
\dot p &=-({\mathcal F}_u+p{\mathcal F}_v)(u,v(u,p),p),\endaligned
$$ \noindent  which amounts to

\begin{equation}\label{eq:selno}
\aligned \dot u &= - \frac{\chi}{2b} u^2 + \frac{\xi}{6b^2}u^3+\frac{c\chi}{2b^2} u^2p+ O(4)\\
\dot p &=-Bu - bp + O(2),\endaligned
\end{equation}

\noindent where   $\xi=12k^2 b^3+ (a_{32}-a_{50}  )b^2+( 3B^2-3BD
-ca_{41})b+ 3 cB(C-k^3).$

 As $\chi = cB-b(C-A +2 k^3) =0$,  the singular point $0$ is  a topological saddle
 of  order three when $b\xi >0$ and is a topological node when $b\xi <0$.

The other singular points of
 the vector field $X_{\mathcal F}=({\mathcal F}_p,p{\mathcal F}_p,-({\mathcal F}_u+p{\mathcal F}_v))$
  are $(0,0,p_1)$ and $(0,0,p_2)$, where $p_i$ are the roots of the equation $bp^2-cp-b=0$,
   which are topological
saddles. In fact,

$$\aligned DX_{\mathcal F}(0,0,p_i)=
\left[\begin{matrix} 0 & -2bp_i+c &0\\ 0 &  p_i(2bp_i-c) & 0 \\
b_1 &b_2 & 3bp_i^2-2cp_i-b\end{matrix}\right],
\endaligned$$
\noindent where $b_{1}=
2Bp_i^2-2p_iC-p_ik^3+p_iA-B-p_i^3k^3+p_i^3C-Dp_i^2$ and
$b_2=p_i^2k^3+2p_i^2C-2p_iD+p_iB-C+k^3+Dp_i^3-p_i^2E$

The non zero eigenvalues
 of $DX_{\mathcal F}(0,0,p_i)$
 are $\lambda_1=-2bp_i^2+cp_i=-(cp_i+2b)=-b(p_i^2+1)$ and
 $\lambda_2=3bp_i^2-2cp_i-b=cp_i+2b=b(p_i^2+1)$.

The implicit surface ${\mathcal F}(u,v,p)=0$ has two  critical
points, one  at $(0,0,p_1)$ and the other at $(0,0,p_2)$. They are
of Morse conic type if and only if  ${\chi}\ne 0$. Otherwise, they
have quadratic co-rank 1 and  Morse index $1$.

Below it will be established  that, at these points, the surface
is locally a topological disk. This will follow from the smooth
right equivalence of ${\mathcal F}$  with $x^2-y^2+z^3=0$,
provided $\xi \neq 0$. The smoothness class ($\geq 1$) is
explained below as a consequence of \cite{me}.

 In fact, at the points $p_1$ and $p_2$, the Hessian of
$\mathcal F$ is given by:

$$\aligned Hess({\mathcal F}(0,0,p_i)=\left[\begin{matrix} 0 & b_{12}&0\\
  b_{12}&-\frac{p_i(cB-Eb+Dc+Ab)}{b}  & c-2bp_i\\
  0 &c-2bp_i & 0\end{matrix}\right]\endaligned$$

\noindent where $ b_{12}=\frac {p_i}b \left( ({k}^{3} -C)c+ (D
-B)b \right).$

 The  eigenvalues of $ Hess({\mathcal F}(0,0,p_i))$ are $0<\mu_1 $,\; $\mu_2<0,$
 with $\mu_1\mu_2=-[b_{12}^2+(c-2bp_i)^2]<0, $ and $\mu_3=0$. The eigenspace
 associated to $\mu_3=0$ is spanned by
  $\ell_3=  \left( ( 2bp_i -c)b, 0,   {p_i} (  ({k}^{3} -C)c+ (D-B)b) \right)$.

 Consider the  implicit equation

 $$({\mathcal F}_v(u,v,p_i+q), {\mathcal F}_p(u,v,p_i+q))=(0,0).$$

 By the  Implicit Function Theorem, the solution of the equation above is a regular curve $c$.
 A regular parametrization of $c$ is given by
 $$c(s)=( b( 2bp_i-c)s, 0(2), p_i+ {p_i} [  ({k}^{3} -C)c+ (D
 -B)b ]s+
 +0(2)).$$

 Therefore, $\gamma(s)={\mathcal F}(c(s))$ is such that $\gamma(0)=p_i$, $\gamma^\prime(0)=0,$
 $\gamma^{\prime\prime}(0)=-p_i b (c^2+4 b^2) \chi  =0$ and
 $\gamma^{\prime\prime\prime}(0)=p_i(2b p_i-c)^3\xi\ne 0.$ This
 follows from the fact that $p_i$ are solutions of $bp^2 -cp + b =0.$

The calculation leading to the expression for
$\gamma^{\prime\prime\prime}(0)$, linear in $\xi$, has been
corroborated by Computer Algebra.

Therefore, by the Generalized Morse Lemma -- see Thom \cite{th},
pp. 75, and  Lopez de Medrano \cite{me}-- there
 exist coordinates $(u_1,v_1,w_1)$ such that
  ${\mathcal F}(u_1,v_1,w_1)=\mu_1 u_1^2+\mu_2 v_1^2+w_1^3$.

Direct examination of proof in this paper gives that the
coordinates are of smoothness class ${k-2-3}\geq 1$, $k-2$ in the
class of $\mathcal F$.

 The phase portrait of $X_{\mathcal F}$ near
the singular point $(0,0,p_i)$ is shown in Fig. \ref{fig:cadeira}.

This picture illustrates the orientation reversing projection  of
the surface into the plane, which occurs in $D^2_{3}$. To get the
orientation preserving case it suffices to look at this picture
upside down.

  \begin{figure}[htbp]
  \begin{center}
 \hskip 1cm
  \includegraphics[angle=0, width=6cm]{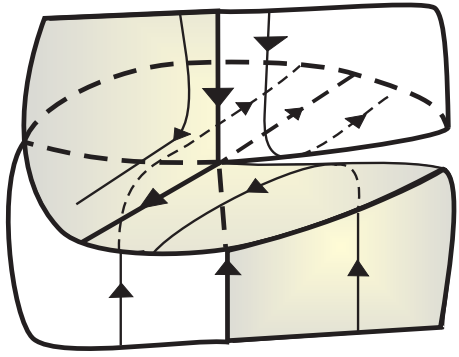}
  \caption{Phase portrait of $X_{\mathcal F} $ in a neighborhood of a
  critical point of corank 1.
  \label{fig:cadeira}}
   \end{center}
  \end{figure}

To conclude the proof it must be noticed that for $X_{\mathcal F}$
the singular point $0$ is a node when $b\xi<0$ and it is a saddle
when $b\xi>0$, which follows from equation (\ref{eq:selno}), while
$p_i$ $(i=1,2)$ are always saddles. This leads to $D^2_{2p}$ in
the first case and to $D^2_{3}$ in the second case. \end{proof}

\subsection  { $D^{2}_{2h}$  Umbilic Points}\label{ssec:d2h}

\begin{proposition}\label{prop:3} Let $p_0$ be an umbilic point and $\alpha(u,v)=(u,v,h(u,v))$ be
  as in equation (\ref{eq:1}).

 If $a=b=0$ and $cB\ne 0$ or if  $b=c=0$ and $aD\ne 0$, then the
principal configuration near $p_0$ is as in  Fig. \ref{fig:d13}.
\end{proposition}

This umbilic point is denoted by $D^{2}_{2h}$.

\begin{proof}
Under these hypotheses, there are two singular points for the
Lie-Cartan vector field \ref{eq:f}; one is located at $p=0$ and
the other at $p= \infty$, that is at $q=0$.

Consider first the implicit differential equation of curvature
lines

\begin{equation}\label{eq:q}
\aligned {\mathcal G}&(u,v,q)=-\frac 12 B u^2 + (k^3-C)uv-\frac 12 Dv^2+O(3)  \\
 +&[  cv+ \frac 12(C+2k^3-A)u^2+(D-B)uv+\frac 12(E-2k^3-C)v^2+O(3)
 ]q\\
 +& [\frac 12 B u^2 -(k^3-C)uv+\frac 12 Dv^2+O(3)]q^2=0.\endaligned
 \end{equation}

\noindent  in the projective chart $(u,v,q)$, where $q=du/dv$.
This implicit surface has a singularity of Morse type at $(0,0,0)$
provided $cB\ne 0$.  In such case, it has the  topological
structure of a cone.

In fact,
$$\aligned Hess({\mathcal G})(0)=
\left[ \begin{matrix} -B & k^3-C& 0\\ k^3-C & -D & c\\ 0 & c & 0
\end{matrix}\right]\endaligned$$ Therefore, $det(Hess({\mathcal
G})(0))=c^2B \ne 0.$
 The Lie-Cartan vector field
 $$Y=(q{\mathcal G}_q,{\mathcal G}_q, -({\mathcal G}_u q+{\mathcal G}_v)) $$
\noindent
 is given by:
$$\aligned
u^\prime &=  q[cv  +\frac 12(C+2k^3-A)u^2+(D-B)uv+\frac 12 (E-2k^3-C)v^2+O(3)]\\
v^\prime &=  cv +\frac 12(C+2k^3-A)u^2+(D-B)uv+\frac 12 (E-2k^3-C)v^2+O(3)\\
q^\prime &=  (C-k^3)u+Dv-cq+\frac 12 cu^2+(c-3k^2c)uv\\
+&\frac 12 a_{14}v^2+(2B-D)uq+(2C-E+k^3)vq+O(3)\endaligned $$

Therefore,
$$\aligned DY(0)= \left[ \begin{matrix}0 & 0 & 0 \\
                    0 & c & 0 \\
                   C-k^3& D & -c \end{matrix}\right] \endaligned $$

\noindent So, the nonzero eigenvalues of $DY(0)$ are $\lambda_1=c$ and $\lambda_2=-c$
and $0$ is a  saddle point for $Y_{\mathcal G}$. The phase portrait of
$Y_{\mathcal G}$ is as shown in Fig. \ref{fig:cone} below.

   \begin{figure}[htbp]
  \begin{center}
   \hskip 1cm
   \includegraphics[angle=0, width=4cm]{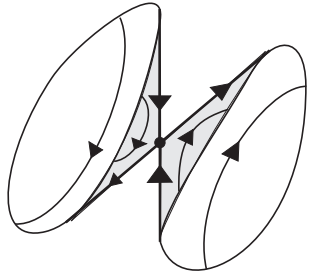}
   \caption{  Cone at $p=\infty $, $q=0$.  \label{fig:cone}}
     \end{center}
   \end{figure}

Consider now the implicit differential equation
\begin{equation}\label{eq:idf} \aligned {\mathcal F}&(u,v,p)=[-\frac 12 B u^2 +
(k^3-C)uv-\frac 12 Dv^2
 -\frac 16a_{41}u^3\\
 -&\frac 12 a_{32}u^2v+\frac{1}2(3k^2c-a_{23})uv^2)-\frac 16 a_{14}v^3+O(4) ]p^2 \\
 +&[  cv+ \frac 12(C+2k^3-A)u^2+(D-B)uv+\frac 12(E-2k^3-C)v^2\\
 +& \frac 16(a_{32}-a_{50})u^3+\frac 12(2k^2c+a_{23}-a_{41})u^2v+\frac 12(a_{14}-a_{32})uv^2\\
 +&\frac
 16(a_{05}-6k^2c-a_{23})v^3+O(4)] p \\
+&  \frac 12 B u^2 -(k^3-C)uv+\frac 12 Dv^2
+\frac 16a_{41}u^3\\
+&\frac 12 a_{32}u^2v+\frac{1}2(a_{23}-k^2c)uv^2+\frac 16
a_{14}v^3+O(4)=0.\endaligned
\end{equation}

\noindent
 in the projective chart $(u,v,p)$, where $p=dv/du$.

 This surface also has a singular point of Morse type at $(0,0,0)$ with
 the topological structure of a cone, provided $ cB\ne 0$.

As above, the Lie-Cartan vector field $Z=({\mathcal F}_p,p{\mathcal F}_p,
 -({\mathcal F}_u +p{\mathcal F}_v)) $ is given by:

\begin{equation}\label{eq:lcp} \aligned
u^\prime &= cv + \frac 12(C+2k^3-A)u^2+(D-B)uv+\frac 12(E-2k^3-C)v^2+O(3).\\
 v^\prime &= p[cv + \frac 12(C+2k^3-A)u^2+(D-B)uv+\frac 12(E-2k^3-C)v^2+O(3)]\\
 p^\prime &=-Bu+(k^3-C)v-cp^2+(A-k^3-2C)up+(B-2D)vp\\
\;&\;\;\;\; - \frac 12a_{41}u^2-a_{32}uv+\frac
12(k^2c-a_{23})v^2+O(3)
\endaligned \end{equation}

In this case the singular point $(0,0,0)$ is not hyperbolic for the
vector field $Z\,=\,Z_{\mathcal F}$. In fact,

$$\aligned DZ(0)=\left[\begin{matrix} 0 & c & 0 \\ 0 & 0 & 0\\ -B & k^3-C & 0
  \end{matrix}\right]\endaligned,$$

\noindent and  the three eigenvalues  are zero.

In order to the describe the phase portrait of $Z$ near zero,
observe that the sheets of the cone, with the origin removed, are
located at $p>0$ and $p<0 $.

To
describe the phase portrait
 of $Z_{\mathcal F}$ near $0$  it is suitable to perform the following
  weighted rescaling of variables.

\begin{equation}\label{eq:wr}
 u=r^2 \bar u, \hskip .5cm v=r^3 \bar v,\hskip .5cm p=r\bar p,
\hskip .5cm \bar u^2+\bar v^2+\bar p^2=1, \hskip .5cm r\geq 0.
\end{equation}

Then, from equation (\ref{eq:idf}) it follows  that

\begin{equation}\label{eq:res} \aligned H(\bar u,\bar v,\bar p)={\mathcal
F}(r^2\bar u, r^3\bar v,r\bar p)= &c\bar v\bar p+\frac 12 B{\bar
u}^2\\
+& r[ (C-k^3)\bar u\bar v+ {\bar u}^2\bar p(k^3+\frac
12C)]+O(r^2)\endaligned
\end{equation}

So, from equations (\ref{eq:lcp}) and (\ref{eq:res}), dropping the
bars  to simplify the notation and dividing by appropriate
 non negative functions,  the following vector field
is obtained: $X=X_u\frac{\partial}{\partial
u}+X_v\frac{\partial}{\partial v}+X_p\frac{\partial}{\partial
p}+X_r\frac{\partial}{\partial r}$, with

$$\aligned u^\prime &= X_u= -2cv^2up+3cv^3+cp^2v+2Bu^2p+2cup^3+rF_1(u,v,p,r)  \\
v^\prime &= X_v= v(2cu^2p+3Bup-3cuv+4cp^3) +r F_2(u,v,p,r) \\
p^\prime &= X_p=  -cuvp-2Bu^3-2cu^2p^2-3Buv^2-4cv^2p^2 +rF_3(u,v,p,r) \\
r^\prime &= X_r =r[ cuv -Bup -c p^3+cv^2p
+rF_4(u,v,p,r)].\endaligned
$$

In the expression above $F_i(u,v,p,0)=0$, ($i=1,\ldots,4$).

The restriction of $X$ to the unitary sphere $u^2+v^2+p^2=1$ is
given by $X_0=X(u,v,p,0)$ which has exactly six singular points,
defined by the intersection of the curves
$$c_1(p)=(0,0,p),\;\;\; c_2(p)=(-\frac{cp^2}B,0,p), \;\;\;
c_3(p)=(-\frac{4cp^2}{3B},-\frac{8cp^3}{9B},p)$$
\noindent with the unitary sphere.

A long, but straightforward, calculation leads to the following
conclusions:
\begin{enumerate}

\item[i)] The non zero eigenvalues of $DX_0(c_1(p))$ are $2cp^3$
and $4cp^3$.

\item[ii)] The non zero eigenvalues of $DX_0(c_2(p))$ are

$-2cp^3 (2c^2p^2+B^2)/B^2$\quad and \quad$cp^3 (2c^2p^2+B^2)/B^2$.

\item[iii)] The non zero eigenvalues of $DX_0(c_3(p))$ are
 $$-\frac
4{81}\frac{ cp^3}{B^2}(27B^2+96c^2p^2+64c^2p^4) \;\;\text{ and}
\;\;-\frac 2{27}\frac{cp^3}{B^2}(27B^2+96c^2p^2+64c^2p^4). $$
\end{enumerate}

Therefore, under the hypothesis $cB\ne 0$, the singular points at
the curves $c_1$ and $c_3$ are hyperbolic (sinks  and sources) and
the singularities at  $c_2$ are hyperbolic saddles.

Also,  the quadratic cone ${\mathcal
Q}_2(u,v,p)=cvp+\frac{B}2u^2=0$ is invariant by $X_0$,  ${\mathcal
Q}_2(c_1(p))={\mathcal Q}_2(c_3(p))=0$ and ${\mathcal
Q}_2(c_2(p))\ne 0$.

Summarizing, the phase portrait of $X_0$ can be described as
follows.

The intersection of the sphere with the cone ${\mathcal
Q}_2^{-1}(0)$ decomposes  the sphere into  two topological disks
$\bar{D}_1\cup \bar{D}_2$ and one cylinder $\bar{C}$.

Each curve  bounding   the disks $\bar{D}_1$ and  $\bar{D}_2$ is
invariant by the flow of $X_0$ and on  it $X_0$ has two hyperbolic
singular points (one attractor and one repeller). Therefore in
each disk the phase portrait of $X_0$ has only one parallel
canonical region.
 In the cylinder $\bar{C}$, $X_0$ has two
hyperbolic saddles and the limit set of ordinary (i.e. non closed)
 orbits, perhaps
with the exception of the saddle separatrices, are the
singularities contained in the bounding curves.

 From  the conclusions above,  proceed  to  the analysis of the
phase portrait of $X$ in the coordinates $(u,v,p,r)$.

 The four singular points on $(c_1(p),0)$ and
$(c_3(p),0)$  $X$ are hyperbolic saddles of stability  index $2$
or $1$; two of each type. In fact the normal derivatives of $X$ in
direction $r$, evaluated at the singularities above, are given by:

$$-cp^3,\;\;\;\; \frac
1{81}\frac{p^3c(27B^2+64p^4c^2+96c^2p^2)}{B^2}.$$

The singular points $(c_2(p),0)$ are not hyperbolic for $X$, but
for the purposes of our analysis it is not necessary to examine
their topological type. The phase portrait of $X$ at the conic
critical point is illustrated in Fig. \ref{fig:bd13cone}.

 \begin{figure}[htbp]
   \begin{center}
  \hskip 1cm
  \includegraphics[angle=0, width=7cm]{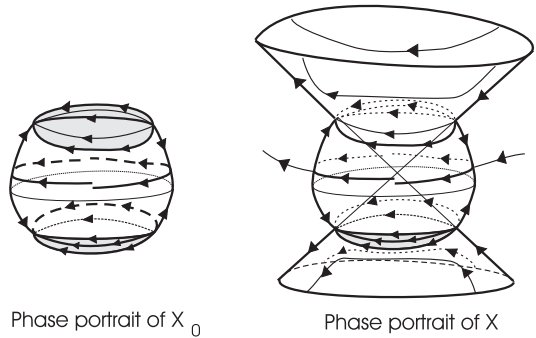}
  \caption{ Blow up at the conic singular point for  $D^{2}_{2h}$  \label{fig:bd13cone}}
   \end{center}
 \end{figure}

Therefore,  near the projective line,  the phase portrait of
$Z\,=\,Z_{\mathcal F}$, restricted to the implicit surface
${\mathcal F}^{-1}(0)$,  is as in Fig. \ref{fig:d13cone}. Here the
point at $\infty$  and that at $0$ have been moved to finite conic
points $p_0$ and $p_{\infty}$, with $p_0 p_{\infty} = -1$, for
better illustration.

 Near the origin,  the projective line is a center manifold with a
quadratic semi-stable singularity. The other separatrix,
asymptotic to the origin,
 sketched in Fig. \ref{fig:d13cone},  also behaves  topologically as a    semi-stable
singularity. This is obtained by blowing down the phase portrait
of $X$  restricted  to  the cone.

\begin{figure}[htpb]
   \begin{center}
  \hskip 1cm
  \includegraphics[angle=0, width=7cm]{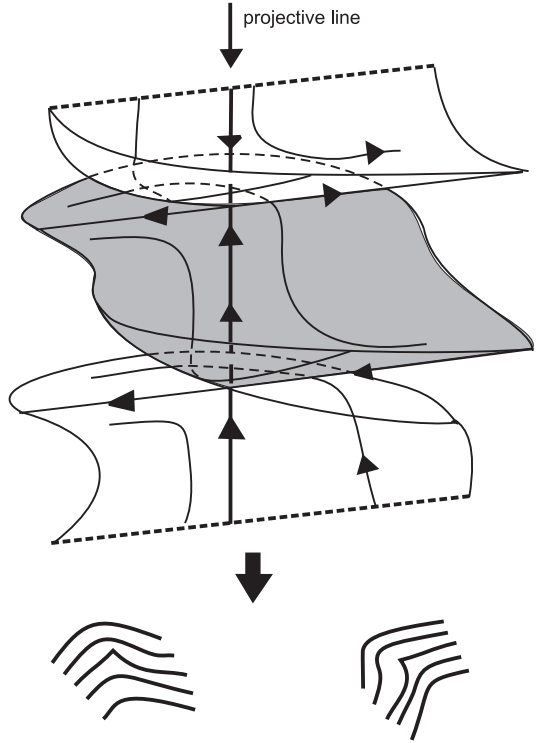}
  \caption{  Lie-Cartan suspension of an umbilic point
  of type $D^{2}_{2h}$  \label{fig:d13cone}}
   \end{center}
 \end{figure}

Gluing together the phase portraits  of $Z\,=\,Z_{\mathcal F}$
near the projective line, in
 the charts $(u,v,p)$ and $(u,v,q)$,  ends the proof.
\end{proof}

\section{A Summarizing  Theorem} \label{sec:4}
 A {\it bi-parametric family of
immersed  surfaces} will be   a family $\alpha_\lambda$ of
immersions of an oriented surface $\mathbb M$ into oriented
$\mathbb R ^3$, indexed by a parameter $\lambda$ ranging in
another surface $\mathbb L$, of {\it parameters}. The families
under consideration will be {\it smooth}  in the sense that
$\alpha(m,\lambda)= \alpha_\lambda (m)$ is  of class $C^\infty$ in
the product manifold ${\mathbb M}\times{\mathbb L}$. The space of
families, denoted $\mathcal F_{\mathbb M \times \mathbb L}$ will
be endowed with the Whitney Fine Topology. See Levine \cite{lt}
and Mather \cite{ma}.

The term {\it generic}, referring to a property of families, means
that it is valid for a collection of families which contains the
intersection of countably many open dense sets in $\mathcal
F_{\mathbb M \times \mathbb L}$. By the appropriate version of
Baire Theorem, these collections are dense. See \cite{lt},
\cite{ma}, where can be found  the basic   terminology  and
concepts that are standard in Singularity Theory. Other reference
on this subject, pertinent to geometric singularities, is Porteous
\cite{po}.

For the sake of simplicity,  from now on, the term {\it smooth}
applied  to mappings and manifolds  will mean  {\it of class}
$C^\infty$.
  See \cite{lt, ma}.

  Consider the space
$\mathbb{J}^{r}(\mathbb{M},\mathbb{R} ^3)$  of $r$-jets of
immersions,
 $\alpha$, of   $\mathbb{M}$  into $\mathbb{R} ^3$, endowed with the
  structure of Principal Fiber Bundle. The base space is  $\mathbb{M}\times \mathbb{R} ^3$; the  fiber
  is the space
  $\mathcal{J}^{r}(2,3)$, where $\mathcal{J}^{r}(2,3)$ is the space of
   $r$-jets of immersions  of $\mathbb{R}^2$ to $\mathbb{R} ^3$,  preserving the respective
     origins. The structure group,  $\mathbb{A}_{+}^{r}$,  is  the product of the group
      $\mathcal{L}_{+}^{r}(2,2)$ of  $r$-jets of origin and  orientation preserving
      diffeomorphisms of $\mathbb{R}^2$, acting on the right by coordinate changes,
       and the group $\mathcal{O}_{+}(3,3)$ of positive
       isometries;
        the action on the left consists  on a  positive rotation of $\mathbb{R}^3$.

Denote by $\Pi_{r,s},\; r\leq s$ the projection of
$\mathcal{J}^{s}(2,3)$ onto $\mathcal{J}^{r}(2,3)$. It is well
known that the group action comutes with projections, \cite{lt}.
For  the present  needs, take $r$ so that $2 \leq r \leq 5$

Recall from Section\ref{sec:2}   that each $5$-jet  of an
immersion at an umbilic point is of the form $(p, P, z)$, with
$(p, P) \in \mathbb{M}\times \mathbb{R} ^3 $ and $z$ is in the
orbit of a polynomial immersion
 $(u,v,h(u,v))$,
where

 \begin{equation}\label{eq:11}
  \aligned h (u,v) &= \frac k2 (u^2+v^2) + \frac a6 u^3 +\frac
b2 u v^2 +\frac c6 v^3
+\frac A{24} u^4 + \frac B6 u^3 v\\
 &+\frac C4 u^2v^2 + \frac D6 u v^3 + \frac E{24} v^2
+\frac{a_{50}}{120}u^5+ \frac{a_{41}}{24}u^4v\\
&+ \frac{a_{32}}{12}u^3v^2 +\frac{a_{23}}{12}u^2v^3+
\frac{a_{14}}{24}uv^4+\frac{a_{05}}{120}v^5
\endaligned
\end{equation}

Recall that the general quadratic part  of $h$ has the form $\frac
{k_{11}}{2} u^2 + k_{12} uv +   \frac {k_{22}}{2} v^2$. Therefore,
the umbilic conditions
 $k_{11}= k_{22}, \, k_{12}=0$, imposed  on  jets, define a closed
submanifold $\mathcal{U}^{5}$, of {\it umbilic jets},  of
codimension $2$ in $\mathcal{J}^{5}(2,3)$.

Recall also that the general cubic part has also a term of the
form $\frac {b^{\prime}}{2} u^2 v$. The expression in equation
(\ref{eq:11}), is a representative of the orbit  the an  umbilic
jets in the space $(a,b,c)$, with   $b^{\prime} =0$.

 \vskip 0.3cm

Define below  the {\it canonic  umbilic stratification} of
$\mathcal{J}^{5}(2,3)$. The term {\it  canonic} means that the
strata are invariant under
 the action of the group
 $\mathbb{A}_{+}^{5}$= $\mathcal{O}_{+}(3,3)\times\mathcal{L}_{+}^{5}(2,2)$.
It is to the orbits of this action that  reference is made below.

\begin{itemize}
\item[1)]
 {\it Umbilic Jets}:  $\mathcal{U}^{5}$, those  in the orbits of  $j^{5}(u,v,h)$,
 where
   $h=h(u,v),$ is as  in equation (\ref{eq:11}), with $k \in \mathbb R$.
It is a closed submanifold of codimension $2$.

\item [2)] {\it Umbilic Jets with Rotationally Symmetric cubic
part}: $\mathcal{Q}^{5}$,  those  in the orbit of $j^{5}(u,v,h)$,
with $h$  in expression (\ref{eq:11}) having $a=b=b^{\prime}=c=0$.
It is a closed submanifold of codimension $4$ inside
$\mathcal{U}^{5}$ and codimension $6$ in $\mathcal{J}^{k}(2,3)$.

\item[3)] {\it Non-Darbouxian Umbilic Jets}: $(\mathcal{ND})^{5}$,
are those umbilic jets outside $\mathcal{Q}^{5}$, which  are in
the orbits of those  whose  expression (\ref{eq:11}) satisfy:

\subitem 3.1) $\tau = b(b-a)= 0$ or \subitem 3.2)
  $b(b-a)\neq 0$ and conditions $D_1$ or $D_2$ in  Theorem \ref{th:1}
  fail.
  This is written
  $\delta_1 = a-2b = 0 $ and
$\delta=c^2 -4b(a - 2b) = 0$.
\end {itemize}

\begin{figure}[htbp]
   \begin{center}
  \hskip 1cm
  \includegraphics[angle=0, width=7cm]{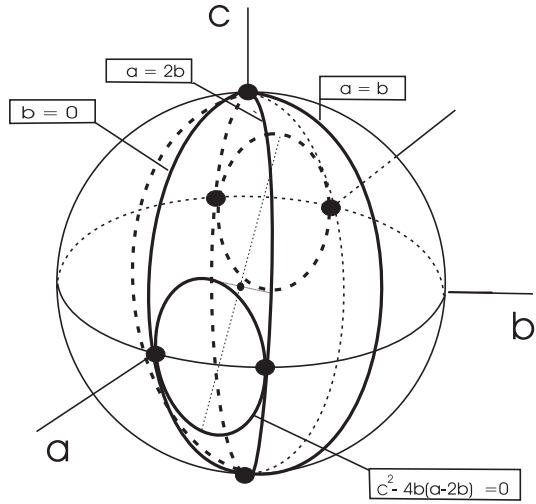}
  \caption{  Stratification of a section of the fiber of $(\mathcal{ND})^{3}$. \label{fig:abc}}
   \end{center}
 \end{figure}

See Fig. \ref{fig:abc} for an illustration of a spherical section
the stratification of the intersection of $(\mathcal{ND})^{3}$,
with the $3$ space $(a,b,c)$, in the notation of equation
(\ref{eq:11}).

The   stratification on this  sphere into points, curves and  open
sets, by the zeros of the homogeneous linear form $\delta_1 =
a-2b$ and quadratic forms $\tau= b(a-b)$ and $\delta=c^2 -4b(a -
2b) =0 $ and their complements is transferred   to
$(\mathcal{ND})^{5}$ by the orbits of the group action, preserving
the codimension of the strata  inside $\mathcal{U}^{5}$. This
codimension will be called  {\it umbilic codimension}, to it  $2$
must be added to get the actual codimension  inside
$\mathcal{J}^{5}(2,3)$.

\begin{itemize}
\item[4)] The variety   $(\mathcal{ND})^{5}$ of {\it
Non-Darbouxian jets}, is further partitioned into

\subitem 4.1)
 $(\mathcal{D}^1_{2,3})^5$, defined by the orbits of jets
with $a=b \neq 0, \; \chi \neq 0$. See Theorem \ref{th:2}.
 \subitem 4.2)
 $(\mathcal{D}^2_{2p})^5$,  defined by the orbits of jets
with $a=b  \neq 0, \; \chi = 0,\;  b\xi < 0$. See Proposition
\ref{prop:2}.

 \subitem 4.3)
   $(\mathcal{D}^2_3)^5$,   defined by the orbits of jets
with  $a=b  \neq 0, \; \chi = 0, \; b\xi > 0$

\subitem 4.4) $(\mathcal{D}^2_{2h})^5$, defined by the orbits of
jets with $a=b = 0, \; cB \neq 0$. This is a manifold of umbilic
codimension $2$. See Proposition \ref{prop:3}.

 \subitem 4.5) $(\mathcal{D}^{1}_{2})^5$ defined by the orbits of jets
 with  $c\tau \neq 0 $ and  $\delta=0$
or  $\delta_1 =0$. See Theorem \ref{th:2}. This is a manifold of
umbilic  codimension $1$.

 \subitem 4.6) $(\mathcal{D}^{2}_{1})^5$ defined by the orbits of jets
 with $c =  0, \, b \neq 0$
and   $\delta_1 =0$.  See Proposition \ref{prop:1}. This is a
manifold of umbilic codimension $2$.
 \subitem 4.7)
$\mathcal{Z}^{5}$,  defined by the orbits of jets with $a=b \neq
0,\; \chi = 0, \;  \xi = 0$, which is a subvariety of umbilic
codimension $3$.

\subitem 4.8) $\mathcal{W}^{5}$,  defined by the orbits of jets
with $a=b \neq 0,\; B=0$, which is a subvariety of umbilic
codimension $3$.
\end {itemize}

\begin{itemize}
\item[5)]The orbits passing through the open sets in  Fig.
\ref{fig:abc} correspond to the Darbouxian umbilic jets:
$(\mathcal{D}_i)^5\, , \,i=1,\, 2,\, 3 $.
\end{itemize}

The {\it canonic stratification} of   $\mathcal{J}^{5}(2,3)$
induces a canonic stratification of
$\mathbb{J}^{5}(\mathbb{M},\mathbb{R} ^3)$ whose strata are
principal sub-bundles, with   codimension equal to that  of their
fibers, which are the canonic strata of $\mathcal{J}^{5}(2,3)$, as
defined   above in items $1, 2, 3, 4 $ and $5$.

The collection of sub-bundles which stratify
$\mathbb{J}^{5}(\mathbb{M},\mathbb{R} ^3)$ will be called  {\it
umbilic stratification}.  The strata are:
$\mathbb{U}^{5}(\mathbb{M},\mathbb{R} ^3)$, corresponding to
$\mathcal{U}^{5}$; $  {\mathbb D}{^5 }_{i}(\mathbb{M},\mathbb{R}
^3), $ $\; i=1, 2, 3$, corresponding to the strata of Darbouxian
umbilic jets $({\mathcal D}_i)^5\, ,\, i=1,\, 2,\, 3 $, and so on,
one bundle for each of the strata in the items above.

\begin{theorem}\label{th:3}
The following properties are generic for smooth bi-parametric
families, $\alpha_\lambda$, of immersions in $\mathcal F_{\mathbb
M \times \mathbb L}$.

 The set $\mathbb U
(\alpha_\lambda)$ to $\mathbb L$ points $(m,\lambda)$ in  $\mathbb
M \times \mathbb L$ such that $m$ is an  umbilic point of
$\alpha_\lambda$ form a smooth regular surface,  which is
stratified as follows:

Darbouxian umbilics,  $D_1$, $D_2$ and $D_3$, occur  on regular
surfaces,

Codimension one umbilics,  $D^{1}_{2}$ and $D^{1}_{2,3}$, occur
along regular curves,

Codimension two  umbilics, $D^{2}_{1}\,$, $D^{2}_{2p} \,$,
$D^{2}_{3}\,$ and $D^{2}_{2h}$, occur at isolated points.

\end{theorem}

\begin{proof}
 Apply   Thom Transversality Theorem in Jet
Spaces, to make  the jet extension $j^{5}\alpha_\lambda$, as a
mapping of $\mathbb M \times \mathbb L$,  transversal to the
canonic stratification of  $\mathbb{J}^{5}(\mathbb{M},\mathbb{R}
^3)$. See Levine \cite{lt} and Mather \cite{ma}.

Notice that the strata  not  listed in the statement of the
Theorem have codimensions larger than $4$ and therefore are
avoided by transversality. While the other strata, listed there,
give the right structure of regular surfaces, curves and points
claimed in the Theorem, since they are obtained by the  pull back
of the canonic strata by  a  $j^{5}\alpha_\lambda$  transversal to
them.

\end{proof}

\begin{remark}\label{rm:fw}
In
 \cite{bif2}, the curves and points above
in Theorem \ref{th:3}  will be  better organized in terms of
Whitney folds and cusps of the projection of  $\mathbb U
(\alpha_\lambda)$ to $\mathbb L$.
 See \cite{lt}.

\end{remark}

\vskip .5cm

\author{\noindent Jorge Sotomayor\\Instituto de Matem\'{a}tica e
Estat\'{\i}stica,\\Universidade de S\~{a}o Paulo, \\Rua do
Mat\~{a}o 1010,
Cidade Universit\'{a}ria, \\CEP 05508-090, S\~{a}o Paulo, S.P., Brazil \\
\\ Ronaldo Garcia\\Instituto de Matem\'{a}tica e
Estat\'{\i}stica,\\Universidade Federal de Goi\'{a}s,\\CEP
74001-970, Caixa Postal 131,\\Goi\^ania, GO, Brazil}

\end{document}